\documentclass{article}[11pt]
\usepackage{amsthm,amsmath,a4wide}
\usepackage{amssymb}
\input xypic
\usepackage[all,tips]{xy}
\usepackage{bbm}

\usepackage%[hidelinks]
[colorlinks,citecolor=blue,linkcolor=black]{hyperref}
\usepackage[inner=1cm,outer=1cm]{geometry}

\newtheorem{thm}{Theorem}[section]
\newtheorem{lem}[thm]{Lemma}
\newtheorem{prop}[thm]{Proposition}
\newtheorem{cor}[thm]{Corollary}
\newtheorem{conj}[thm]{Conjecture}

\theoremstyle{definition}
\newtheorem{defin}[thm]{Definition}

\theoremstyle{remark}
\newtheorem{remark}[thm]{Remark}
\newtheorem{example}[thm]{Example}

\newcommand{\bth}{\begin{thm}}
\renewcommand{\eth}{\end{thm}}
\newcommand{\bpr}{\begin{prop}}
\newcommand{\epr}{\end{prop}}
\newcommand{\ble}{\begin{lem}}
\newcommand{\ele}{\end{lem}}
\newcommand{\bco}{\begin{cor}}
\newcommand{\eco}{\end{cor}}
\newcommand{\bde}{\begin{defin}}
\newcommand{\ede}{\end{defin}}
\newcommand{\bex}{\begin{example}}
\newcommand{\eex}{\end{example}}
\newcommand{\bre}{\begin{remark}}
\newcommand{\ere}{\end{remark}}
\newcommand{\bcj}{\begin{conj}}
\newcommand{\ecj}{\end{conj}}

\newcommand{\beq}{\begin{equation}}
\newcommand{\eeq}{\end{equation}}

\newcommand{\ot}{{\otimes}}
\newcommand{\op}{{\oplus}}
\newcommand{\lb}{\label}
\newcommand{\nl}{\newline}
\newcommand{\bpf}{\begin{proof}}
\newcommand{\epf}{\end{proof}}

\newcommand{\da}{{\text-}}

\newcommand{\E}{{\cal E}}
\newcommand{\A}{{\cal A}}

\newcommand{\C}{{\cal C}}
\newcommand{\Z}{{\cal Z}}
\newcommand{\D}{{\cal D}}
\newcommand{\B}{{\cal B}}
\newcommand{\V}{{\cal V}}

\newcommand{\R}{{\cal R}}
\renewcommand{\S}{{\cal S}}
\newcommand{\I}{{\cal I}}

\newcommand{\M}{{\cal M}}
\newcommand{\N}{{\cal N}}

\newcommand{\F}{{\cal F}}
\newcommand{\Vect}{{\cal V}{\it ect}}

\begin{document}
\author{Alexei Davydov}
\title{Tensor categories}
\maketitle
\date{}

\begin{center}\small
Department of Mathematics, Ohio University, Athens, OH 45701, USA \\  davydov@ohio.edu 
\end{center}

\begin{abstract}
We present here definitions and constructions basic for the theory of monoidal and tensor categories. 
We provide references to the original sources, whenever possible. 
Group-theoretical categories are used as examples. 
\end{abstract}
\bigskip

{\bf Keywords/Objectives:} monoidal categories and functors, braided and symmetric monoidal categories, rigid and spherical monoidal categories, dimension, algebras and modules in monoidal categories, linear and abelian monoidal categories, tensor and fusion  categories.

\tableofcontents

\section{Introduction}

Monoidal categories were defined in \cite{be} and further studied in \cite{mc0} and  \cite{ke}. 
The notion turned out to be a fundamental one. 
From philosophical point of view it is the result of applying the so-called {\em categorification} (in the sense of e.g. \cite{cf,bd}) to the notion of a monoid, i.e. the process of replacing sets with categories, functions with functors, and equations between functions by natural isomorphisms between functors, which in turn should satisfy certain equations of their own, called ``coherence laws".  The pseudo-monoid nature of monoidal categories explains their appearance in low dimensional topology, through (extended) topological field theories \cite{tu}.
Other (somewhat related) applications of monoidal categories include mathematical logic and theoretical computer science (see \cite{bs}). 

A special kind of monoidal categories is formed by linear (hom-sets are vector spaces, composition and tensor product of morphisms are bilinear maps) and abelian categories. The natural source of such categories is the representation theory of groups and other related structures (Lie algebras, Hopf algebras, etc.). 
An abelian monoidal  (or {\em tensor}) category can thus be seen as the category of all possible linear realisations of a certain symmetry. This point of view is beneficial in such areas as the theory of motives or Langlands correspondence and was driving the development of (symmetric) tensor categories \cite{dm, de1}. 

The concept of {\em fusion} (finite semi-simple tensor) categories originated in high energy physics \cite{ms}.
Braided fusion categories appear as representation categories of chiral algebras in 2d conformal field theory or as values on the circle of corresponding 3d topological field theories. More general fusion categories are formed by boundary data of such theories. 
Currently the largest consumer of fusion categories is condensed matter physics (topological states of matter) \cite{ki}. 

For definitions and basic facts of category theory, such as adjoint functors, limits, ends see \cite{mc}.
We follow Mac Lane's notation denoting the set of morphisms between objects $X$ and $Y$ of a category $\A$ by $\A(X,Y)$. 
For a functor $F:\C\to\D$ we denote by  $\I m_f(F)$ the full subcategory of $\D$ with objects of the form $F(X)$ for $X\in\C$ (the {\em full image} of $F$) and 
by $\I m(F)$ is the subcategory of $\D$ with objects of the form $F(X)$ for $X\in\C$ and morphisms of  the form $F(f)$ for $f$ a morphism in $\C$ (the {\em proper image} of $F$). 
Thus $F$ factors as $\C \to \I m(F) \to \I m_f(F) \to\D$, where the first functor is full, the second is faithful, and the last is fully faithful (an {\em embedding}).

\section{Monoidal categories}

A category $\C$ is {\em semi-groupal} if it comes equipped with a functor (the {\em tensor product})
$$\C\times\C\to \C\qquad (X,Y)\mapsto X\ot Y$$
together with a collection of isomorphisms (the {\em associativity constraint}) $a_{X,Y,Z}:X\ot(Y\ot Z)\to (X\ot Y)\ot Z$ natural in $X,Y,Z\in\C$, 
such that the diagram
\beq\lb{ac}\xygraph{
!{<0cm,0cm>;<3.7cm,0cm>:<0cm,2.5cm>::}
!{(0,.6)}*+{(X\ot Y)\ot(Z\ot W)}="t"
!{(-1,0)}*+{X\ot (Y\ot(Z\ot W))}="l"  !{(1,0)}*+{((X\ot Y)\ot Z)\ot W}="r"
!{(-.6,-.8)}*+{X\ot ((Y\ot Z)\ot W)}="lb"  !{(.6,-.8)}*+{(X\ot (Y\ot Z))\ot W}="rb"
%%%%%%%
"l":"t" ^{a_{X,Y,Z\ot W}}
"t":"r" ^{a_{X\ot Y,Z,W}}
"l":"lb" _{1\ot a_{Y,Z,W}}
"lb":"rb" ^{a_{X,Y\ot Z,W}}
"rb":"r" _{a_{X,Y,Z}\ot 1}
}\eeq
commutes for any $X,Y,Z,W\in\C$. 
\nl
A semi-groupal category is {\em strict} if $a_{X,Y,Z}=1$ for all $X,Y,Z\in\C$. 
\bex
For a category $\M$ denote by $\E nd(\M)$ the category of endofunctors $\M\to\M$. This is a strict semi-groupal category with composition of functors as the tensor product.  
\eex
A functor $F:\C\to\D$ between semi-groupal categories is {\em semi-groupal} if it comes equipped with a collection of isomorphisms (the {\em semi-groupal constraint})  $F_{X,Y}:F(X)\ot F(Y)\to F(X\ot Y)$ natural in $X,Y\in\C$, 
such that the diagram
\beq\lb{fc}\xygraph{
!{<0cm,0cm>;<4.6cm,0cm>:<0cm,1.6cm>::}
!{(-1,1)}*+{F(X)\ot(F(Y)\ot F(Z))}="tl"  !{(0,1)}*+{F(X)\ot F(Y\ot Z)}="t"  !{(1,1)}*+{F(X\ot (Y\ot Z))}="tr"
!{(-1,0)}*+{(F(X)\ot F(Y))\ot F(Z)}="bl"  !{(0,0)}*+{F(X\ot Y)\ot F(Z)}="b"  !{(1,0)}*+{F((X\ot Y)\ot Z)}="br"
%%%%%%%
"tl":"bl" _{a_{F(X),F(Y),F(Z)}}
"tl":"t" ^(.55){1\ot F_{Y,Z}}
"t":"tr" ^{F_{X,Y\ot Z}}
"bl":"b" ^(.55){F_{X,Y}\ot 1}
"b":"br" ^{F_{X\ot Y,Z}}
"tr":"br" ^{F(a_{X,Y,Z})}
}\eeq
commutes for any $X,Y,Z\in\C$. 
\nl
A natural transformation $c:F\to G$ between semi-groupal functors $F,G:\C\to\D$ is {\em semi-groupal} if the diagram
\beq\lb{mn}\xygraph{
!{<0cm,0cm>;<2.6cm,0cm>:<0cm,1.4cm>::}
!{(-1,1)}*+{F(X)\ot F(Y)}="tl"  !{(0,1)}*+{F(X\ot Y)}="tr"  
!{(-1,0)}*+{G(X)\ot G(Y)}="bl"  !{(0,0)}*+{G(X\ot Y)}="br"  
%%%%%%%
"tl":"bl" _{c_X\ot c_Y}
"tl":"tr" ^(.55){F_{X,Y}}
"bl":"br" ^(.55){G_{X,Y}}
"tr":"br" ^{c_{X\ot Y}}
}\eeq
commutes for any $X,Y\in\C$. 
\nl
Clearly composites of semi-groupal functors and natural transformations are semi-groupal.
\nl
%More precisely, the above definitions are of {\em semi-groupal} category, functor, and natural transformation correspondingly. 
A semi-groupal category $\C$ is {\em monoidal} if it has an object $I\in \C$ (the {\em monoidal unit}) together with collections of isomorphisms $l_X:I\ot X\to X$, $r_X:X\ot I\to X$ natural in $X\in\C$ such that the diagram
\beq\lb{uc}\xygraph{
!{<0cm,0cm>;<2cm,0cm>:<0cm,1cm>::}
!{(-1,1)}*+{X\ot(I\ot Y)}="tl"  !{(1,1)}*+{(X\ot I)\ot Y}="tr"  
!{(0,0)}*+{X\ot Y}="b"  
%%%%%%%
"tl":"b" _{1\ot l_Y}
"tl":"tr" ^{a_{X,I,Y}}
"tr":"b" ^{r_X\ot 1}
}\eeq
commutes for any $X,Y\in\C$. 
\bre
The monoidal unit object is unique up to an isomorphism, i.e. being monoidal is a property of a semi-groupal category. 
\nl
In a strict monoidal category one also has $l=r=1$. 
\ere
A semi-groupal functor $F:\C\to\D$ is {\em monoidal} if it comes equipped with an isomorphism $\iota:I\to F(I)$ such that the diagrams
\beq\lb{uf}
\xygraph{
!{<0cm,0cm>;<2.6cm,0cm>:<0cm,1.3cm>::}
!{(-1,1)}*+{F(I)\ot F(X)}="tl"  !{(0,1)}*+{F(I\ot X)}="tr"  
!{(-1,0)}*+{I\ot F(X)}="bl"  !{(0,0)}*+{F(X)}="br" 
%%%%%%%
"bl":"tl" ^{\iota\ot 1}
"tl":"tr" ^(.55){F_{I,X}}
"tr":"br" ^{F(l_X)}
"bl":"br" ^{l_{F(X)}}
}\qquad\qquad
\xygraph{
!{<0cm,0cm>;<2.6cm,0cm>:<0cm,1.3cm>::}
!{(-1,1)}*+{F(X)\ot F(I)}="tl"  !{(0,1)}*+{F(X\ot I)}="tr"  
!{(-1,0)}*+{F(X)\ot I}="bl"  !{(0,0)}*+{F(X)}="br" 
%%%%%%%
"bl":"tl" ^{1\ot \iota}
"tl":"tr" ^(.55){F_{X,I}}
"tr":"br" ^{F(r_X)}
"bl":"br" ^{r_{F(X)}}
}\eeq
commute for any $X\in\C$. 

The {\em scalars} of a monoidal $\C$ form the submonoid $\C_s(I,I)=\{c\in\C(I,I)|\ l_X(c\ot 1_X) = r_X(1_X\ot c)\ \forall X\in\C\}$ of the monoid of endomorphisms $\C(I,I)$. By the Eckman-Hilton argument the monoid $\C(I,I)$ is commutative. 

For a monoidal $\C$ the {\em opposite category} $\C^{op}$ (the category with morphism sets $\C^{op}(X,Y) = \C(Y,X)$) is monoidal with $a^{op}_{X,Y,Z} = a_{X,Y,Z}^{-1},\ l^{op}_X = l_X^{-1},\ r^{op}_X = r_X^{-1}$. 
The {\em reverse} $\C^{rev}$ of a monoidal $\C$ is the category $\C$ with a new tensor product $X\ot^{rev}Y = Y\ot X$ and $a^{rev}_{X,Y,Z} = a_{Z,Y,X}^{-1},\ l^{rev}_X = r_X,\ r^{rev}_X = l_X$. 
We also write: $\C^{oprev}=(C^{rev})^{op}$. 

\section{Braided and symmetric monoidal categories}

A monoidal category $\C$ is {\em braided} if it comes equipped with a collection of isomorphisms
$c_{X,Y}:X\ot Y\to Y\ot X$ ({\em braidings}) natural in $X,Y\in\C$ such that the diagrams
\beq\lb{bc}
\xygraph{ !{0;/r3.5pc/:;/u3.5pc/::}[]*+{X\ot(Y\ot Z)} 
(
  :[u(.7)r(.5)]*+{(X\ot Y)\ot Z} ^(.4){a_{X,Y,Z}}
  :[r(2)]*+{Z\ot(X\ot Y)} ^{c_{X\ot Y,Z}}
  :[d(.7)r(.5)]*+{(Z\ot X)\ot Y}="r" ^(.6){a_{Z,X,Y}}
,
  :[d(.7)r(.5)]*+{X\ot(Z\ot Y)} _(.4){1c_{Y,Z}}
  :[r(2)]*+{(X\ot Z)\ot Y} ^{a_{X,Z,Y}^{-1}}
  : "r" _(.6){c_{X,Z}1}
)}
\qquad\qquad
\xygraph{ !{0;/r3.5pc/:;/u3.5pc/::}[]*+{(X\ot Y)\ot Z} 
(
  :[u(.7)r(.5)]*+{X\ot (Y\ot Z)} ^(.4){a_{X,Y,Z}^{-1}}
  :[r(2)]*+{(Y\ot Z)\ot X} ^{c_{X,Y\ot Z}}
  :[d(.7)r(.5)]*+{Y\ot(Z\ot X)}="r" ^(.6){a_{Y,Z,X}^{-1}}
,
  :[d(.7)r(.5)]*+{(Y\ot X)\ot Z} _(.4){c_{X,Y}1}
  :[r(2)]*+{Y\ot (X\ot Z)} ^{a_{Y,X,Z}^{-1}}
  : "r" _(.6){1 c_{X,Z}}
)
}\eeq
 commute for all $X,Y,Z\in\C$. It is also assumed that $c_{I,Y} = l_Y^{-1}r_Y$ and $c_{X,I} = r_X^{-1}l_X$. 
 \nl
A monoidal functor $F:\C\to\D$ between braided monoidal categories is {\em braided} if the diagram
\beq\lb{cbf}\xymatrix{F(X)\ot F(Y) \ar[rr]^{F_{X,Y}} \ar[d]_{c_{F(X),F(Y)}} && F(X\ot Y) \ar[d]^{F(c_{X,Y})} \\
G(X)\ot G(Y) \ar[rr]^{G_{X,Y}} && G(X\ot Y)}\eeq
commutes for all $X,Y\in\C$ \cite{js1}. 
\nl
The {\em conjugate} $\overline\C$ of a braided monoidal category $\C$ is the monoidal category $\C$ with the braiding
$\overline{c}_{X,Y}=c_{Y,X}^{-1}$. 
\bex
For a monoidal category $\C$ define $\Z(\C)$ as the category of pairs $(Z,z)$, where $Z\in\C$ together with a collection of isomorphisms $z_{X}:X\ot Z\to Z\ot X$ natural in $X\in\C$ ({\em half-braiding}) and such that the diagram obtained by replacing $c$ with $z$ in the first diagram of \eqref{bc} commutes. 
The category $\Z(\C)$ is monoidal with the tensor product $(Z,z)\ot (W,w) = (Z\ot W,z|w)$, where $(z|w)_X:X\ot(Z\ot W)\to (Z\ot W)\ot X$ is defined by a diagram similar to the second diagram of \eqref{bc}. 
The monoidal category $\Z(\C)$ is braided (with the braiding $c_{(Z,z),(W,w)} = z_W$)
and is known as the {\em monoidal} (or {\em Drinfeld}) {\em centre} of $\C$ \cite{js,ma}.
\eex
A braided monoidal category $\C$ is {\em symmetric} if $c_{Y,X}c_{X,Y} = 1$ for all $X,Y\in\C$.
The full subcategory $\Z_{sym}(\C) = \{X\in\C|\ c_{Y,X}c_{X,Y} = 1\ \forall Y\in\C\}$ is symmetric (the {\em symmetric or Mueger's} centre of $\C$). 

A braiding of $\C$ can be thought of as the structure of a monoidal functor $\C\to\C^{rev}$ on the identity functor $Id_\C$. 
For a symmetric $\C$ the composite $\C\to\C^{rev}\to \C$ coincides with $Id_\C$ as a monoidal functor.
\nl
More generally, a {\em balancing} of a braided monoidal $\C$ is a monoidal isomorphism between the composite $\C\to\C^{rev}\to \C$ and $Id_\C$ \cite{js1}. The group $Aut_\ot(Id_\C)$ of monoidal automorphisms of the identity functor acts freely on the set of possible balancings of a balanced category. 
\nl
Not all structures of a monoidal functor $\C\to\C^{rev}$ on the identity functor $Id_\C$ come from braidings of $\C$.
A monoidal $\C$ together with a structure of monoidal functor $\C\to\C^{rev}$ on the identity functor $Id_\C$ such that the composite $\C\to\C^{rev}\to \C$ coincides with $Id_\C$ as a monoidal functor was called {\em coboundary} in \cite{dr}. 

\section{Rigid monoidal categories}

The {\em (left) internal hom} $[X,Y]$ of objects $X,Y$ of a semi-groupal category $\C$ is the terminal object $([X,Y], ev)$ in the category of pairs $(T,\tau)$ where $T$ is an object and $\tau:T\ot X\to Y$ is a morphism in $\C$ correspondingly. 
Equivalently, $[X,-]$ is the right adjoint to the functor $-\ot X$, i.e. there is a collection of isomorphisms
$\C(T,[X,Y])\simeq \C(T\ot X,Y)$ natural in $X,Y,T\in\C$. 
The universal property of the internal hom gives rise to  morphisms
$$[Y,Z]\ot[X,Y]\longrightarrow [X,Z],\qquad X\ot[Y,Z] \stackrel{}{\longrightarrow} [Y,X\ot Z],\qquad [X,[Y,Z]] \stackrel{\sim}{\longrightarrow} [X\ot Y,Z]$$
natural in $X,Y,Z\in\C$. 
A semi-groupal category $\C$ is {\em (left) closed} if $[X,Y]$ exists for any $X,Y\in\C$ \cite{ek}. 
\bre
The existence of $[X,-]$, as a right adjoint, can e.g. be deduced from the adjoint functor theorems (see \cite{mc0}). 
If $[X,-]$ exists the functor $-\ot X$ is a left adjoint and hence preserves colimits. Thus the tensor product of a closed category preserves colimits (in both arguments).
\ere
In a monoidal category $\C$ one has $[I,X]\simeq X$ for any $X\in\C$. 
A monoidal category $\C$ is {\em (left) rigid} of {\em (left) autonomous} if it is (left) closed and $[I,Y]\ot[X,I]\longrightarrow [X,Y]$ is an isomorphism for any $X,Y\in\C$. In this case $X^*=[X,I]$ is called the {\em (left) dual} of $X$. 
The (left) dual $X^*$ of $X$ is characterised by the existence of morphisms $ev:X^*\ot X\to I, coev:I\to X\ot X^*$ such that $(coev\ot 1)(1\ot ev) = 1_X$ and $(1\ot coev)(ev\ot 1) = 1_{X^*}$ (here we suppress associators). Such morphisms are unique up to an automorphism of $I$ (see \cite{sr}). 
\bex
The (left) dual to an endofunctor $F$ as an object of the monoidal category $\E nd(\M)$ is its left adjoint.
Thus a left rigid $F\in\E nd(\M)$ preserves limits and a right rigid endofunctor preserves colimits in $\M$. 
\eex
\bre
A monoidal functor $F:\C\to \D$ preserves duals $F(X^*)=F(X)^*$. 
Thus for a rigid $\C$ the proper and full images $\I m(F)$ and $\I m_f(F)$ of a monoidal functor  are rigid. 
\ere
\bre
The monoidal functor $\C\to\E nd(\C)$ sending $X$ to the endofunctor $Y\mapsto X\ot Y$ sends rigid object to functors preserving limits and colimits. Thus the tensor product of a rigid category preserves limits and colimits (in both arguments).
\ere

For rigid $\C$ the assignment $X\mapsto X^*$ extends to morphisms and defines a monoidal equivalence $\C\to\C^{oprev}$. 
A {\em pivotal structure} on a rigid monoidal $\C$ is a monoidal isomorphism $\tau$ between $Id_\C$ and the composite $(-)^{**}:\C\to\C^{oprev}\to \C$ \cite{fy}. 
The group $Aut_\ot(Id_\C)$ acts freely on the set of pivotal structures of a pivotal category. 
\nl
For an endomorphism $f:X\to X$ of a pivotal category $\C$ the composite $tr(f)=ev_{X^*}(\tau_X f\ot 1)coev_X$ is an endomorphism of $I$ called the {\em trace} of $f$ (here we suppress associators).
A pivotal category $\C$ is {\em spherical} if $tr(f^*) = tr(f)$ for any endomorphism $f$ in $\C$ \cite{bw}.

\section{Module categories}

A category $\M$ is a {\em (left) module} category over a monoidal category $\C$ if it comes equipped with a monoidal functor $\C\to\E nd(\M)$. Equivalently, a $\C$-module structure of $\M$ is encoded in a functor 
$\C\times\M\to \M,\ (X,M)\mapsto X*M$. 
The monoidal constraint of $\C\to\E nd(\M)$ corresponds to a collection of isomorphisms (the {\em action associativity constraint}) $a_{X,Y,M}:X*(Y* M)\to (X\ot Y)* M$ natural in $X,Y\in\C$, $M\in\M$, 
such that the diagrams similar to \eqref{ac},\eqref{uc} commute. 
A functor $F:\M\to\N$ between $\C$-module categories is {\em $\C$-module functor} if it comes equipped with a collection of isomorphisms (the {\em module constraint})  $F_{X,M}:F(X)* F(M)\to F(X* M)$ natural in $X\in\C$, $M\in\M$, 
such that the diagrams similar to \eqref{fc},\eqref{uf} commute. 
A natural transformation $c:F\to G$ between $\C$-module functors $F,G:\M\to\N$ is {\em $\C$-module} if the diagram similar to \eqref{mn} commutes.
\nl
Clearly composites of $\C$-module functors and natural transformations are $\C$-module.
In particular, $\C$-module endofunctors of a $\C$-module category $\M$ together with $\C$-module natural transformations form a strict monoidal category $\E nd_{\C\da}(\M)$.  

\bex
The tensor product can be seen as an action of a monoidal category $\C$ on itself (the {\em regular} left $\C$-module category). 
The assignment $F\mapsto F(I)$ is a monoidal equivalence $\E nd_{\C\da}(\C)\to\C$. Thus any monoidal category is equivalent to a strict one (Mac Lane's coherence theorem). 
\eex

Similarly one defines right $\C$-module and $\C$-bimodule categories. 
For a $\C$-bimodule category $\M$ denote by $\E nd_{\C\da\C}(\M)$ the strict monoidal category of $\C$-bimodule endofunctors. 

\bex
The tensor product turns a monoidal category $\C$ into a bimodule over itself. 
The assignment $F\mapsto F(I)$ extends to a monoidal equivalence $\E nd_{\C\da\C}(\C)\to\Z(\C)$ with the $\C$-bimodule structure of $F$ corresponding to the half braiding on $F(I)$. 
\eex

The notion of internal hom extends to module categories.
The {\em internal action hom} $[M,N]_\M$ of objects $M,N$ of a (left) $\C$-module category $\M$ is the terminal object $([M,N]_\M, ev)$ in the category of pairs $(T,\tau)$ where $T\in\C$ and $\tau:T* M\to N$ is a morphism in $\M$. 
Equivalently, $[M,-]_\M$ is the right adjoint to the functor $-*M$, i.e. there is a collection of isomorphisms
$\C(T,[M,N]_\M)\simeq \M(T*M,N)$ natural in $M,N\in\M,T\in\C$ \cite{os}. 
The universal property of the internal hom gives rise to  morphisms
$$[M,N]_\M\ot[L,M]_\M\longrightarrow [L,N]_\M,\qquad X\ot[M,N]_\M \stackrel{}{\longrightarrow} [M,X*N]_\M,\qquad [X,[M,N]_\M] \stackrel{\sim}{\longrightarrow} [X*M,N]_\M$$
natural in $L,M,N\in\M, X\in\C$. 
For rigid $\C$ the second morphism is an isomorphism. 

\section{Internal algebra}

The language of algebras and their modules internal to a monoidal category is a useful tool for defining and computing with module categories and for constructing braided monoidal categories.

An {\em associative unital algebra} in a monoidal category $\C$ is a triple $(A,\mu,\iota)$ consisting of an object $A\in\C$ together with a {\em multiplication} $\mu:A\otimes A\to A$ and a {\em unit} map $\iota:I\to A$, satisfying the {\em associativity} $\mu\left(\mu\otimes 1_A\right) = \mu\left(1_A\otimes\mu\right),$ and the {\em unit} $\mu\left(\iota\otimes 1_A\right) = 1_A = \mu\left(1_A\otimes\iota\right)$ axioms. %Where it will not cause confusion,  we will be talking about an algebra $A$, suppressing its multiplication and unit maps.\newline
\nl
A {\em right module} over an algebra $A$ is a pair $(M,\nu)$, where $M$ is an object of $\C$ and $\nu:M\otimes A\to M$ is a morphism ({\em action map}), such that $\nu\left(\nu\otimes 1_A\right) = \nu\left(1_M\otimes\mu\right).$ A {\em homomorphism} of right $A$-modules $M\to N$ is a morphism $f:M\to N$ in $\C$ such that $\nu_N\left(f\otimes 1_A\right) = f\nu_M.$
\nl
Right modules over an algebra $A\in\C$ together with module homomorphisms form a category $\C_A$. The forgetful functor $\C_A\to\C$ has a left adjoint, which sends an object $X\in\C$ into the {\em free $A$-module} $X\otimes A$, with $1_X\ot \mu$ as the $A$-module structure. %defined by $\xymatrix{X\otimes A\otimes A \ar[r]^{1\mu} & X\otimes A.}$
%Since the action map $M\otimes A\to M$ is an epimorphism of right $A$-modules, any right $A$-module is a quotient of a free module.\newline
More generally, for $X\in\C$ and $M\in\C_A$ the tensor product $X\ot M$ is a right $A$-module under $1_X\ot \nu_M$.
This categorical pairing is a left $C$-module action on $\C_A$. The free $A$-module functor $-\ot A:\C\to\C_A$ and the forgetful functor $\C_A\to\C$ are $C$-module functors. 
\nl
Conversely, let $M$ be an object of a $\C$-module category $\M$ such that internal action homs $[M,-]_\M$ exist, 
Then $[M,M]_\M$ is an algebra and $[M,N]_\M$ is a right $[M,M]_\M$-module  in $\C$.
For rigid $\C$ the functor $[M,-]_\M:\M\to \C_{[M,M]_\M}$ is $\C$-module and the free $[M,M]_\M$-module functor factorises
$$\xymatrix{\C \ar[rr] \ar[rd]_{-*M} && \C_{[M,M]_\M} \\ & \M\ar[ru]_{[M,-]_\M} }$$

An algebra $A$ in a braided monoidal category is {\em commutative} if $\mu c_{A,A}=\mu$. 
For $M,N\in\C_A$ denote by $M\ot_AN$ the coequaliser of $\nu_M1,(1\nu_N)(1c_{A,N}):M\ot A\ot N\to M\ot N$. 
When these coequalisers exist for all $M,N\in\C_A$ the category $\C_A$ becomes monoidal with the tensor product $\ot_A$ \cite{pa}. The free-module functor $\C\to\C_A$ is monoidal and {\em central}, i.e. lifts to a braided monoidal functor $\C\to\Z(\C_A)$.
\nl
Conversely for a central functor $F:\C\to\D$ with the right adjoint $R:\D\to\C$ the object $R(I)\in\C$ has a structure of commutative algebra. The functor $F$ induces a monoidal functor $\tilde F:\C_{R(I)}\to\D$ and $F = \tilde F\circ(-\ot R(I))$ \cite{dmno}. Note that the algebra $R(I)$ coincides with the internal action hom $[I,I]$ for the $\C$-action on $\D$ defined by the functor $F$. 
\bex
For a $\C$-module category $\M$ there is defined a central monoidal functor $\Z(C)\to\E nd_\C(\M)$ (called {\em $\alpha$-induction}). The corresponding internal action hom $[Id_\M,Id_\M]$ is a commutative algebra $Z(\M)$ in $\Z(C)$ called the {\em full centre} of $\M$.
For $\M=\C_A$ the full centre $Z(A) = Z(\C_A)$ is the terminal object of the category of the pairs $(Z,\zeta)$, where $Z\in\Z(\C)$ and $\zeta:Z\to A$ is a morphism in $\C$ such that $\mu(\zeta\ot 1_A) = \mu(1_A\ot\zeta)z_A$ (where $z_A$ is the half-braiding of $Z$) \cite{ffrs,da}.
\eex

More generally for a commutative algebra $A\in\Z(\C)$ the coequaliser $M\ot_AN$ (with $c_{A,N}$ replaced by the half braiding of $A$) is a tensor product on $\C_A$ \cite{sh}.

For closed $\C$ the category of $A$-modules $\C_A$ is also closed with the internal hom $[M,N]_A$ being by the equaliser of $[1,\widetilde\nu_N],[\nu_M,1]:[M,N]\to [M,[A,N]]$, where $\widetilde\nu_N:N\to[A,N]$ is induced by the module structure $\nu_N:N\ot A\to N$ \cite{pa}. 
For a rigid $\C$ the category $\C_A$ is also rigid if the algebra $A$ is {\em separable}, i.e. an algebra $A$ such that the {\em canonical (trace) pairing} $ev_A z_{A^*}(\mu\ot 1)(1\ot coev_A)\mu:A\ot A\to A$ is non-degenerate \cite{ffrs}.

For a commutative algebra $A$ in a braided monoidal category $\C$ the monoidal category $\C_A$ is not braided in general. However it contains a distinguished full braided subcategory $\C_A^{loc}$ formed by so-called local modules.
An $A$-module $(M,\nu)$ is {\em local} if $\nu c_{A,M}c_{M,A} =\nu$ \cite{pa,ffrs}. 
Note that the full embedding $\C_A^{loc}\to \C_A$ is a central monoidal functor.

%Hopf algebras

\section{Linear and abelian monoidal categories}

It is convenient when dealing with additive and linear categories to use the language of categories {\em enriched} in a symmetric monoidal category $\V$, or {\em $\V$-category}. 
In a $\V$-category $\A$ homs $\A(X,Y)$ are not sets but rather objects of $\V$ and composition is a morphism $\A(Y,Z)\ot_\V\A(X,Y)\to\A(X,Z)$ (see \cite{ke1} for all details). 
For $\V$ being the category of sets $\S et$, $\V$-category is just a category.
When $\V$ is the category of abelian groups, $\V$-category is an additive category ($\A(X,Y)$ is an abelian group and composition is bilinear).
When $\V$ is the category $\Vect$ of vector spaces over a fixed field $k$, $\V$-category is a linear category ($\A(X,Y)$ is a vector space and composition is bilinear).

A $\V$-category $\C$ is {\em monoidal} if it comes equipped with a $\V$-functor $\C\boxtimes_\V\C\to \C$ together with an associativity constraint and the unit object. Here $\A\boxtimes_\V\B$ is the $\V$-category with objects $Ob(\A\boxtimes_\V\B) = Ob(\A)\times Ob(\B)$ and homs $(\A\boxtimes_\V\B)((X,Z),(Y,W)) = \A(X,Y)\ot_\V\B(Z,W)$.
For example in a monoidal linear category the tensor product on morphisms is given by bilinear maps $\C(X,Y)\ot_k\C(Z,W)\to\C(X\ot Z,Y\ot W)$. 
The addition operation on morphisms can be quite useful. 
For an additive functor $F:\C\to\D$ the proper image $\I m(F)$ can be described as the quotient $\C/Ker(F)$ by the {\em ideal} $Ker(F)$ of all morphisms mapped by $F$ to zero. When the functor $F$ is monoidal the ideal $Ker(F)$ is also {\em monoidal}, i.e. $Ker(F)$ is closed not just under compositions but also under tensor products with arbitrary morphisms. 
\bex
In a spherical linear category negligible morphisms form a monoidal ideal $N$. A morphism $f:X\to Y$ is {\em negligible} if 
$tr(fg)=0$ for any $g:Y\to X$ \cite{tu,bw}. 
\eex

A {\em pre-abelian} category is an additive category with finite limits and colimits ({\em finitely bicomplete}).
A pre-abelian category is {\em abelian} if any monomorphism is a kernel and any epimorphism is a cokernel. 

For (co)complete $\V$ (with (co)continuous tensor product) the category $\widehat\A = \F unct_\V(\A^{op},\V)$ of $\V$-functors (the {\em $\V$-presheaf category}) is also (co)complete. 
The {\em Yoneda embedding} $Y:\A\to \widehat\A$ sending $X$ to the functor $\A(X,-)$ is fully faithful.
\nl
For a monoidal $\V$-category $\C$ the presheaf category $\widehat\C$ is monoidal with respect to the 
{\em Day convolution} defined as the coend $(U\ot V)(X) = \int^{Y,Z\in\C}\C(X,Y\ot Z)\ot_\V U(Y)\ot_\V V(Z)$
\cite{dy}.
If $\V$ is closed the monoidal presheaf category $\widehat\C$ is also closed with the internal hom
given by the end $[U,V](X) = \int_{Y\in\C}[U(Y),V(X\ot Y)]_\V$. 
\bex
For a linear category $\A$ the presheaf category $\widehat\A$ is abelian.

\eex
An abelian linear monoidal category $\C$ is {\em tensor} if it is rigid and the identity object $I$ is simple and $\C(I,I)=k$. 
Here an object is {\em simple} if all morphisms out of it are monomorphisms and morphism in it are epimorphisms. 

\section{Finite tensor and fusion categories}

A tensor category is {\em finite} if all hom-spaces are finite dimensional and any object has a finite length (i.e. has a finite filtration with simple factors). As an abelian category a finite tensor category is equivalent to the category of finite dimensional modules over a finite dimensional algebra. 
As the result a finite tensor category is finitely complete and cocomplete, and a tensor functor between finite tensor categories has left and right adjoints. 
In particular, internal action homs for a finite module category exist. 

The {\em categorical dimension} $dim(\C)$ of a finite tensor category is the dimension $d(Z(I))$ of the full centre of the trivial algebra $I\in\C$. 
The categorical dimension has the following properties:
$dim(Z(\C)) = dim(\C)^2$, $dim(\C_A)d(A) = dim(\C)$ and $dim(\C_A^{loc})d(A)^2 = dim(\C)$ for a separable commutative algebra $A\in\Z(\C)$ such that $\C(I,A) = k$  \cite{ffrs,egno}. 

A finite tensor category $\C$ is {\em incompressible} if any tensor functor $\C\to\D$ into a tensor category is fully faithful (see \cite{ceo} and references therein). 
For an incompressible finite tensor category $\C$ there should be no non-trivial separable commutative algebras $A$ in $\Z(\C)$ with $\C(I,A) = k$.
%Braided tensor categories have their version of {\em incompressiblity} with tensor functors $\C\to\D$ assumed to be braided 
\bex
An incompressible symmetric tensor category in zero characteristic is either the category of vector spaces $\Vect$ or its super-analog $s\Vect$ \cite{de1}.
A symmetric tensor category $\C$ is {\em Tannakian} if there is a faithful symmetric tensor functor $\C\to\V ect$ and {\em super-Tannakian} if there is a faithful symmetric tensor functor $\C\to s\V ect$. 
\nl
The structure of incompressible symmetric tensor categories in finite characteristic is much richer (see \cite{ceo} and references therein).
\eex
For an incompressible braided finite tensor category $\C$ there should be no non-trivial separable commutative algebras $A$ in $\C$ with $\C(I,A) = k$ ({\em complete anisotropy} of $\C$, see \cite{dmno}). 
\bex
A braided tensor category $\C$ is {\em non-degenerate} if $\Z_{sym}(\C) = \Vect$ and is {\em slightly-degenerate} if $\Z_{sym}(\C) = s\Vect$.
\nl
Completely anisotropic braided finite tensor categories are either non-degenerate or slightly-degenerate.
\eex

A finite tensor category is {\em fusion} if it is semi-simple, i.e. any object is a direct sum of simple objects. 
\bex
Equivalence classes of completely anisotropic non-degenerate braided fusion  categories form a group, the so-called {\em Witt group} \cite{dmno}.
The same is true for slightly-degenerate braided fusion  categories. 
\eex

\section{Example: pointed categories}

Let $G$ be a semi-group. 
Denote by $\V(G)$ the category of vector spaces (over a field $k$) graded by $G$ with grading preserving linear maps as morphisms. 
Note that it is a semi-simple category (any object is a direct product of simple objects) with simple objects being one dimensional vector space $I(g)=k$ concentrated in the a single degree $g\in G$.
It is also a monoidal category with the tensor product of vector spaces graded by $(U\ot_k V)_f = \op_{gh=f}U_g\ot_k V_h$. 
The unit object $I$ is the one dimensional vector space concentrated in the trivial degree $I=I(e)$. 
It follows from the naturality that the general form of an associativity constraint on homogeneous $u\in U_f, v\in V_g, w\in W_h$ is $a_{U,V,W}(u\ot (v\ot w)) = \alpha(f,g,h)(u\ot v)\ot w$, where $\alpha(f,g,h)\in k^\times$ are non-zero scalars. 
The coherence \eqref{ac} is equivalent to the function $\alpha:G^{\times 3}\to k^\times$ being a 3-cocycle of $G$ with coefficients in the trivial $G$-module $k^\times$. 
Denote by $\V(G,\alpha)$ the category of $G$-vector spaces with the associator corresponding to a 3-cocycle $\alpha\in Z^3(G,k^\times)$. 
\nl
The category $\V(G,\alpha)$ is closed with the left internal hom $[U,V]=Hom_k(U,V)$ being the space of linear maps with the grading $[U,V]_f=\{l|\ l(U_g)\subset V_{fg}\}$ (the right internal hom $\{U,V\}=Hom_k(U,V)$ is graded by $\{U,V\}_f=\{l|\ l(U_g)\subset V_{gf}\}$). A graded vector space $U$ is dualisable if and only if $U$ is finite dimensional and supported on invertible elements of $G$.
For a (finite) group $G$ the subcategory $\V_{fd}(G,\alpha)$ of finite dimensional graded vector spaces is rigid and fusion (a so-called {\em pointed fusion category}). 
\nl
A monoidal functor $F:\V(G,\alpha)\to \V(Q,\beta)$ corresponds to a group homomorphism $\varphi:G\to Q$ (the effect of $F$ on simple objects) and a 2-cochain $\phi\in C^2(G,k^\times)$ (defining the monoidal constraint). Explicitly,  $F(U)_q=\op_{\varphi(g)=q}U_g$ with the monoidal constraint $F_{U,V}(u\ot v) = \phi(f,g)u\ot v$ on $u\in U_f, v\in V_g$.
The coherence \eqref{fc} is equivalent to the coboundary condition $d(\phi) = \alpha^{-1}\varphi^*(\beta)$ in $C^*(G,k^\times)$. 

Let $X$ be a $G$-set with the action $G\times X\to X,\ (g,x)\mapsto g.x$. 
The category $\M(X)$ of $X$-graded vector spaces is a module category over $\V(G)$ with the categorical action $U*M=U\ot_kM$ graded by $(U*M)_x=\op_{g.y=x}U_g\ot_kM_y$ and the identity associativity constraint. 
To make the category of $X$-graded vector spaces a module category over $\V(G,\alpha)$ one needs a coboundary $\xi\in C^2(G,Map(X,k^\times))$ for $\alpha$ considered an element of $Z^3(G,k^\times)\subset Z^3(G,Map(X,k^\times))$. Here $Map(X,k^\times)$ is the multiplicative group of maps $X\to k^\times$ with the $G$-action coming from the $G$-set structure. Such $\xi$ defines an action associativity constraint $U*(V*M)\to (U\ot V)*M,\ \ u\ot(v\ot m)\mapsto \xi(f,g,x)(u\ot v)\ot m$ with $u\in U_f, v\in V_g, m\in M_x$. The action coherence is equivalent to the equation $d(\xi)=\alpha$. 
The  $\V(G,\alpha)$-module category $\M(X,\xi)$ is action closed with the action internal hom $[M,N]=\op_g[M,N]_g$  where  $[M,N]_g=\{l\in Hom_k(M,N)|\ l(U_x)\subset V_{g.x}\}$. 
Note that any semi-simple module category over $\V(G,\alpha)$ is (equivalent to) such $\M(X,\xi)$ for some $(X,\xi)$ (with $X$ being the set of isomorphism classes of simple objects). 
Indecomposable semi-simple $\V(G,\alpha)$-module categories correspond to transitive $G$-sets and have the form $\M(G/H,\xi)$ for a subgroup $H\subset G$. Shapiro's lemma allows one to replace $\xi\in C^2(G,Map(G/H,k^\times))$ by $\tilde\xi\in C^2(H,k^\times)$ with $d(\tilde\xi)=\alpha|_H$. 
Categories of module endofunctors $\E nd_{\V(G,\alpha)}(\M(G/H,\xi))$ are called {\em group-theoretical} (see \cite{egno} and references therein). For example, $\E nd_{\V(G,1)}(\M(G/G,1))$ is the category $\R ep(G)$ of representations of $G$. 
\nl
An associative algebra in $\V(G,\alpha)$ is a $G$-graded vector space $A=\op_gA_g$ with a multiplication, which is {\em graded} (or $\alpha$-){\em associative}, i.e. $a(bc) = \alpha(f,g,h)(ab)c$ for any $a\in A_f, b\in A_g, c\in A_h$. The subalgebra $A_e=\V(G,\alpha)(I,A)$ is an associative algebra in $\V ect$. 
For example, the action internal end $[I(x),I(x)]$ for a simple object $I(x)\in \M(X,\xi)$ is supported on the stabiliser $H=St_G(x)$ and coincides with the twisted group algebra $k[H,\tilde\xi]$. The category of modules $\V(G,\alpha)_{k[H,\tilde\xi]}$ is equivalent to $\M(G/H,\xi)$ as a $\V(G,\alpha)$-module category. 
\nl
The monoidal centre $\Z(\V(G,1))$ is equivalent to the category with objects, $G$-graded vector spaces $Z=\op_gZ_g$ with a linear $G$-action compatible via $f(Z_g) = Z_{fgf^{-1}}$. The tensor product is the standard one.  
The half-braiding of $Z$ with a $G$-graded vector space $V$ is given by $z_{U}(u\ot z) = f.z\ot u$ for $u\in U_f$.
The monoidal centre $\Z(\V(G,\alpha))$ has a similar model (see e.g. \cite{ds} for details). 

The category $\V(G,\alpha)$ affords a braiding only if the group $G$ is abelian.
The naturality restricts the general form of a braiding to $c_{U,V}(u\ot v) = \gamma(f,g)v\ot u$ for $u\in U_f, v\in V_g$, with some $\gamma(f,g)\in k^\times$. 
The coherence \eqref{bc} is equivalent to $(\alpha,\gamma)$ being an {\em abelian} 3-cocycle (a 3-cocycle of the Eilenberg-Mac Lane cohomology complex $C^*(G,2,k^\times)$). 
The coherence \eqref{cbf} for the monoidal constraint of a braded monoidal functor $F:\V(G,\alpha,\gamma)\to \V(Q,\beta,\xi)$ is equivalent to the coboundary condition $d(\phi) = (\alpha,\gamma)^{-1}\varphi^*(\beta,\xi)$ in $C^*(G,2,k^\times)$. 
In particular,  up to a braided equivalence $\V(G,\alpha,\gamma)$ depends only on the cohomology class $[(\alpha,\gamma)]\in H^3(G,2,k^\times)$. Eilenberg and Mac Lane proved that the assignment $(\alpha,\gamma)\mapsto q(g)=\gamma(g,g)$ gives an isomorphism between $H^3(G,2,k^\times)$ and the group $Q(G,k^\times)$ of {\em quadratic} functions, i.e. functions $q:G\to k^\times$ such that $q(g^{-1})=q(g)$ and $q(f)q(g)q(h)q(fgh) = q(fg)q(fh)q(gh)$ (see e.g. \cite{js,dmno} for details and references).  It is convenient to denote (the braided equivalence class of)  $\V(G,\alpha,\gamma)$ by $\C(G,q)$. For example, the  conjugate $\overline{\C(G,q)}$ is $\C(G,q^{-1})$ and the symmetric centre $\Z_{sym}(\C(G,q)) = \C(Ker(q),1)$, where $Ker(q)=\{g\in G|\ q(fg)=q(f)\ \forall f\in G\}$ is the {\em kernel} of the quadratic function $q$. The braided linear category $\C(G,q)$ is non-degenerate whenever the quadratic function $q$ is.
\nl
The twisted group algebra $k[H,\tilde\xi]$ in $\C(G,q)$ is commutative if and only if $H$ is an {\em isotropic} subgroup, i.e. $q|_H=1$. Since $\tilde\xi$ is uniquely determined by $H$ we denote the commutative algebra $k[H,\tilde\xi]$ by $A(H)$. 
The category $\C(G,q)_{A(H)}^{loc}$ of local $A(H)$-modules is equivalent to $\C(H^\perp/H,q)$, where $H^\perp=\{g\in G|\ q(fg)=q(f)\ \forall f\in H\}$ is the {\em orthogonal complement} of $H$. 
The braided linear category $\C(G,q)$ is completely anisotropic whenever the pair $G,q)$ is, i.e. whenever $G$ does not have proper isotropic subgroups.
The subgroup of the Witt group generated by braided pointed fusion categories is the classical Witt group of quadratic functions studied by C.T.C. Wall) see \cite{dmno} for details).

\section{Conclusions}

The development of the theory of tensor categories was driven by the applications in
representation theory (resulting in a good understanding of Tannakian categories \cite{dm}),
high energy and condensed matter physics (bringing fusion categories to the centre of attension).
Concepts crucial for the emergent theory of tensor categories came from or play an important role in these applications:
non-degenerate braided fusion  categories (representation categories of chiral algebras), module categories (categories of boundary data), full centres (full state spaces of conformal field theories, modular invariants), Witt equivalence (boson condensation, defects in 2d conformal field theories).
\nl
Interest to non-rational conformal field theories in two dimensions stimulates the study of more general (finite) tensor categories.
\nl
Higher categorical analogues of tensor categories (fusion 2-categories) play an important role in 4d topological field theories.
Progress in their study is impossible without better understanding of their partial decategorifications. 
Thus it is natural to expect further developments in the area of (finite set-theoretical) monoidal categories.

Due to the space limitation many important references are not listed here.

{\footnotesize

\bibliographystyle{ams-alpha}

}
\end{document}